\documentclass[12pt]{amsart}
\usepackage{amsmath,amscd,amssymb,amsfonts}
\newcommand{\h}{\hbox}
\newcommand{\q}{\quad}
\newcommand{\qt}{\,\,\,}
\newcommand{\qf}{\,\,\,\,}
\newcommand{\nin}{\noindent}
\newcommand{\bs}{\par\bigskip}
\newcommand{\ms}{\par\medskip}
\newcommand{\sk}{\par\smallskip}
\newcommand{\msum}{\h{$\sum$}}
\newcommand{\mopl}{\h{$\bigoplus$}}
\newcommand{\ssb}{\raise.15ex\h{${\scriptscriptstyle\bullet}$}}
\newcommand{\ssc}{\,\raise.15ex\h{${\scriptstyle\circ}$}\,}
\newcommand{\kd}{\hbox{$\frac{k}{d}$}}
\newcommand{\nt}{\hbox{$\frac{n}{2}$}}
\newcommand{\ot}{\hbox{$\frac{1}{2}$}}
\newcommand{\tauY}{\tau_{\,Y}}
\newcommand{\C}{{\mathbf C}}
\newcommand{\N}{{\mathbf N}}
\newcommand{\PP}{{\mathbf P}}
\newcommand{\Q}{{\mathbf Q}}

\newcommand{\Z}{{\mathbf Z}}
\newcommand{\D}{{\mathcal D}}
\newcommand{\Hc}{{\mathcal H}}
\newcommand{\Ic}{{\mathcal I}}
\newcommand{\OO}{{\mathcal O}}
\newcommand{\X}{{\mathcal X}}
\newcommand{\XX}{{}\,\overline{\!{\mathcal X}}{}}
\newcommand{\f}{{}\,\overline{\!f}{}}
\newcommand{\al}{{\alpha}}
\newcommand{\alt}{\widetilde{\alpha}}
\newcommand{\Ct}{\widetilde{C}}
\newcommand{\Ht}{\widetilde{H}}

\newcommand{\Kfb}{K_f^{\ssb}}
\newcommand{\sN}{{}^s\!N}
\newcommand{\dd}{{\partial}}
\newcommand{\ddd}{{\rm d}}
\newcommand{\df}{{\rm d}f}
\newcommand{\mm}{{\mathfrak m}}
\newcommand{\ga}{{\gamma}}
\newcommand{\la}{\lambda}
\newcommand{\Om}{\Omega}
\newcommand{\Si}{\Sigma}
\newcommand{\Sg}{{\rm Sing}}
\newcommand{\Sp}{{\rm Sp}}
\newcommand{\Gr}{{\rm Gr}}
\newcommand{\bl}{\bigl}
\newcommand{\br}{\bigr}
\newcommand{\into}{\hookrightarrow}
\newcommand{\simto}{\buildrel\sim\over\longrightarrow}
\newcommand{\ges}{\geqslant}
\newcommand{\les}{\leqslant}
\newcommand{\ODP}{\hbox{\rm ODP}}
\newcommand{\Pn}{{}^P\!n}
\begin{document}
\title[Theorems of Griffiths and Steenbrink]
{Generalization of theorems of Griffiths\\
and Steenbrink to hypersurfaces\\
with ordinary double points}
\author[A. Dimca]{Alexandru Dimca}
\address{Universit\'e C\^ote d'Azur, CNRS, LJAD, Nice, France.}
\email{dimca@unice.fr}
\author[M. Saito]{Morihiko Saito}
\address{RIMS Kyoto University, Kyoto 606-8502 Japan}
\email{msaito@kurims.kyoto-u.ac.jp}
\begin{abstract}
Let Y be a hypersurface in projective space having only ordinary double points as singularities. We prove a variant of a conjecture of L. Wotzlaw on an algebraic description of the graded quotients of the Hodge filtration on the top cohomology of the complement of Y except for certain degrees of the graded quotients, as well as its extension to the Milnor cohomology of a defining polynomial of Y for degrees a little bit lower than the middle. These partially generalize theorems of Griffiths and Steenbrink in the Y smooth case, and enable us to determine the structure of the pole order spectral sequence. We then get quite simple formulas for the Steenbrink and pole order spectra in this case, which cannot be extended even to the simple singularity case easily.
\end{abstract}
\maketitle
\centerline{\bf Introduction}
\bs\nin
Let $Y$ be a hypersurface in $X=\PP^n$. Consider the following condition:
$$\h{$Y$ has only ordinary double points as singularities.}
\leqno(\ODP)$$
Let $R:=\C[x_0,\dots,x_n]$ with $x_0,\dots,x_n$ the coordinates of $\C^{n+1}$. Let $f\in R$ be a defining polynomial of $Y$. Set $U=X\setminus Y$, and
$$d:=\deg f,\q m:=\bl[\nt\br],\q J:=(\dd f/\dd x_0,\dots,\dd f/\dd x_n)\subset R,\q I:=\sqrt{J}\subset R.$$
Here $J$ is called the Jacobian ideal of $f$, and $I$ is the graded ideal consisting of finite sums of homogeneous polynomials vanishing at $\Sg\,Y\subset X$ if condition~$(\ODP)$ is assumed.
Let $R_k$ denote the degree $k$ part of $R$, and similarly for
$I_k$, etc. We have the following.
\ms\nin
{\bf Conjecture 1} (L. Wotzlaw \cite[6.5]{Wo}). Under the assumption $(\ODP)$ we have
$$\Gr_F^pH^n(U,\C)=(I^{q-m+1}/I^{q-m}J)_{(q+1)d-n-1}\q(q=n-p\in\Z).$$
\sk
Here $F$ is the Hodge filtration as in \cite{De}, and $I^j=R$ for $j\les 0$.
This is a generalization of Griffiths' theorem on rational integrals \cite{Gri} in the $Y$ smooth case (see also \cite{Gre}). The following is known:
\ms\nin
{\bf Theorem~1} (\cite[Theorem~2.2]{DSW}). {\it Conjecture~$1$ holds if $q\les m$, that is, if $p\ges n-m$.}
\ms
This theorem was actually proved in the case of general singularities by modifying $m$ and $I$ appropriately. (More precisely, $m$ is replaced by $\alt_Y$ in (1.1.3), see also (1.1.6).)
\sk
Let $\Kfb:=(\Om^{\ssb},\df\wedge)$ be the Koszul complex associated with the action of $\df\wedge$ on the algebraic differential forms $\Om^{\ssb}:=\Gamma(\C^{n+1},\Om_{\C^{n+1}}^{\ssb})$.
It is a graded complex with $\deg x_i=\deg \ddd x_i=1$, and $\df\wedge$ is a morphism of degree $d$. Set
$$\sN:=H^n(\Kfb),\q M:=H^{n+1}(\Kfb),\q M':=H_{\mm}^0M,\q M'':=M/M',
\leqno(0.1)$$
where $\mm\subset R$ is the maximal ideal generated by the $x_i$, and $H_{\mm}^0$ is the local cohomology. These are graded $R$-modules.
In \cite{DS2} we used $N=\sN(-d)$ instead of $\sN=N(d)$ with $(d)$ a shift of grading.
Under the assumption $(\ODP)$ we have the isomorphisms (see for instance \cite{Di3}):
$$M_k=(R/J)_{k-n-1},\q M''_k=(R/\sqrt{J})_{k-n-1}=(R/I)_{k-n-1}\q(k\in\Z),
\leqno(0.2)$$
\sk
Conjecture~1 is naturally extended to the case of the Milnor cohomology $H^n(f^{-1}(1),\C)$, generalizing Steenbrink's theorem \cite{St1} in the isolated singularity case, at least for lower degrees $q$. There is a technical difficulty as is explained in \cite[Section 1.8]{DS1} if one tries to generalize directly the argument in the proof. However, this can be avoided by using the Thom-Sebastiani type theorem, and we get the following (see (2.4) below).
\ms\nin
{\bf Theorem~2.} {\it Under the assumption $(\ODP)$ the pole order spectral sequence in $(1.2.4)$ below degenerates at $E_2$, and there are canonical isomorphisms
$$\Gr_F^pH^n(f^{-1}(1),\C)_{{\bf e}(-k/d)}=M_k=(R/J)_{k-n-1}\qf\h{for}\qf p=\bl[n+1-\kd\br],\,\,\kd\les\nt,
\leqno(0.3)$$
where ${\bf e}(-k/d):=\exp(-2\pi ik/d)$, and $H^n(f^{-1}(1),\C)_{\la}$ is the $\la$-eigenspace under the action of the monodromy. Moreover the Hodge filtration $F$ on the left-hand side can be replaced with the pole order filtration $P$.}
\ms
Indeed, the assertion with $F$ replaced by $P$ follows from Theorem~(2.1) below, and we can show the coincidence of $F$ and $P$ in the case of Theorem~2 by using the Thom-Sebastiani theorems, see (2.4) below. The isomorphism (0.3) for $\kd\les 1$ is already known by the relation with the multiplier ideals, see \cite{mult}.
Without assuming condition~$(\ODP)$, the isomorphism in (0.3) holds for $\kd\les\alt_Y$ with $\alt_Y$ as in (1.1.3), see Remark~(2.10) below.
\sk
Now consider the Steenbrink spectrum (\cite{St2}, \cite{St3}) and the pole order spectrum \cite{DS2}:
$$\Sp(f)=\msum_{\al}\,n_{f,\al}t^{\al},\q\Sp_P(f)=\msum_{\al}\,\Pn_{f,\al}t^{\al},$$
see (1.4) below for the definition. Let $\ga_k$ be the integers defined by
$$\msum_k\,\ga_k\,t^k=(t+\cdots+t^{d-1})^{n+1}.$$
Note that $\msum_k\,\ga_k\,t^{k/d}$ coincides with the Steenbrink spectrum of a homogeneous polynomial with an isolated singularity, see \cite{St1}. Since the Euler characteristic of a finite dimensional complex is independent of the differential, we have
$$\dim M_k-\dim\sN_{k-d}=\ga_k.
\leqno(0.4)$$
\sk
By using Theorem~(2.1) below together with (0.4), we get the following (see (2.5) below).
\ms\nin
{\bf Theorem~3.} {\it Set $^s\nu_k:=\dim\sN_k\,\,(k\in\Z)$. Under the assumption~$(\ODP)$ we have}
$$\Sp_P(f)=\bl(\,t^{1/d}+\cdots+t^{(d-1)/d}\,\br)^{n+1}-\sum_{nd/2<k\les nd/2+d}{}^s\nu_k\,t^{k/d}-\sum_{k>nd/2+d}(^s\nu_k-{}^s\nu_{k-d})\,t^{k/d}.$$
\ms
Here $nd$ must be even in the case $\kd=\nt$ or $\nt+1$. Note that
$${}^s\nu_k-{}^s\nu_{k-d}\in[0,\tauY]\q\h{with}\q\tauY:=\#|\Sg\,Y|,$$
since $\bl\{{}^s\nu_k\br\}_k$ is a weakly increasing sequence with values in $[0,\tauY]$ (by using (0.6) below). The assertion for $\kd=\nt$ or $\nt+1$ is closely related to the following.
\ms\nin
{\bf Proposition~1} (\cite[Chapter 6, Theorem~4.5]{Di1}). {\it Under the assumption $(\ODP)$ we have}
$$\dim H^{n-1}(f^{-1}(1),\C)=\begin{cases}{}^s\nu_{nd/2}&\h{if}\qt nd\qt\h{is even},\\ \,0&\h{if}\qt nd\qt\h{is odd}.\end{cases}$$
\sk
This also follows from Theorem~(2.1) below. Indeed, it is essentially equivalent to the vanishing of the morphism $\ddd^{(1)}:\sN_{nd/2}\to M_{nd/2}$ induced by the differential $\ddd$.
\sk
As for the Steenbrink spectrum, we have a quite simple formula as follows (see (2.6) below).
\ms\nin
{\bf Theorem~4.} {\it Set $k_0:=[nd/2]$. Under the assumption~$(\ODP)$ we have}
$$\Sp(f)=\bl(\,t^{1/d}+\cdots+t^{(d-1)/d}\,\br)^{n+1}-\tauY\,t^{k_0/d}\,\bl(\,t^{1/d}+\cdots+t^{d/d}\,\br).$$
\ms
This follows from Theorem~2 for the coefficients $n_{f,\alpha}$ of $\Sp(f)$ with $\alpha=\kd\les\nt$. We have a partial symmetry of the Steenbrink spectrum by using a spectral sequence associated with the weight filtration on the vanishing cycle sheaf (1.5), and this implies the assertion for $\alpha=\kd>\nt+1$. For the remaining case we calculate the Euler characteristic of $\PP^n\setminus Y$, see (2.6) below. Theorem~4 is compatible with a formula for the spectrum of a hyperplane arrangement \cite[Theorem~3]{BS} in the case $n=2$.
\sk
By Theorems~3 and 4, we get the following relation between the Steenbrink and pole order spectra.
\ms\nin
{\bf Corollary~1.} {\it Let $p(k)\in\Z$ with $\kd-p(k)\in\bl(\nt,\nt+1\br]$. Under the assumption~$(\ODP)$ we have}
$$\Sp(f)-\Sp_P(f)=\msum_{k/d>n/2+1}\,\bl({}^s\nu_k-{}^s\nu_{k-d}\br)\,\bl(t^{k/d}-t^{k/d-p(k)}\br).$$
\ms
Indeed, this immediately follows from Theorems~3 and 4, since
$$\msum_{p>0}\,\bl({}^s\nu_{k+pd}-{}^s\nu_{k+pd-d}\br)=\tauY-{}^s\nu_k.$$
However, Corollary~1 does not give a formula for $\dim\Gr_F^p\Gr_P^{p+q}H^n(f^{-1}(1),\C)_{\la}$ by combining it with Proposition~(2.9) and (1.4.3) below, which imply a formula for the relation between the two spectra. Indeed, there may be cancellations on the right-hand side of (1.4.3).
\sk
Using the improved version of \cite[Theorem~5.3]{DS2} explained in \cite[Remark 5.6(i)]{DS2} together with Proposition~1, we also get the following (see (2.7) below).
\ms\nin
{\bf Theorem~5.} {\it Let $M'_k$ be as in $(0.1)$. Under the assumption $(\ODP)$ there is a canonical injective morphism
$$M'_k\into\Gr_P^pH^n(f^{-1}(1),\C)_{{\bf e}(-k/d)}\q\h{for}\qt p=\bl[n+1-\kd\br],\,\,k\in\Z_{>0},$$
which is induced by the canonical surjection $M_k\to M_k^{(\infty)}$ in the notation of $(1.2.7)$ below together with the action of $\dd_t^{n-p}$ in $(1.2.3)$ and the isomorphism $(1.2.5)$.}
\ms
This seems to be related with \cite{B1}, \cite{B2}. Theorem~5 and the first assertion of Theorem~2 can be extended to the case $Y$ has only weighted homogeneous isolated singularities, see \cite{wh}.
It is unclear whether Theorem~5 holds with the pole order filtration $P$ replaced by the Hodge filtration $F$ (except for the case $k/d\les n/2$ by Theorem~2). However, combining Theorem~4 with results of \cite{DiSt1}, \cite{Di4} (see also Theorem~9 below), we get the following (see (2.8) below).
\ms\nin
{\bf Proposition~2.} {\it Let $M'_k$ be as in $(0.1)$. Under the assumption $(\ODP)$ there is an inequality
$$\dim M'_k\les\dim\Gr_F^pH^n(f^{-1}(1),\C)_{{\bf e}(-k/d)}\q\h{for}\qt p=\bl[n+1-\kd\br],\,\,k\in\Z_{>0},$$
where the equality holds if $\kd\in\bl(\nt,\nt+1\br]$.}
\ms
In this paper we also treat a variant of Wotzlaw conjecture studied in \cite{DSW}.
Let $\Ic\subset\OO_X$ denote the reduced ideal of $\Sg\,Y\subset X$. Set
$$I^{(i)}_k:=\Gamma(X,\Ic^i(k)),\q I^{(i)}=\mopl_k\,I^{(i)}_k.$$
We have the inclusions $(I^i)_k\subset I^{(i)}_k$ together with the equalities $(I^i)_k=I^{(i)}_k$ for $k\gg 0$ although these equalities do not always hold in general, see \cite[Section 2.3]{DSW}.
By definition we have exact sequences
$$0\longrightarrow I^{(i)}_k\longrightarrow R_k\buildrel{\beta_{k}^{(i)}}\over\longrightarrow\mopl_{y\in\Sg\,Y}\,\OO_{X,y}/\mm_{X,y}^i,
\leqno(0.5)$$
choosing a trivialization of $\OO_{X,y}(k)$, where $\mm_{X,y}$ is the maximal ideal of $\OO_{X,y}$. We have a variant of Conjecture~1 as follows:
\ms\nin
{\bf Conjecture 2.} Under the assumption $(\ODP)$ we have
$$\Gr_F^pH^n(U,\C)=(I^{(q-m+1)}/I^{(q-m)}J)_{(q+1)d-n-1}\q(q=n-p\in\Z).$$
\sk
Note that Theorem~1 implies Conjecture~2 for $q\les m$.
In \cite{DSW}, the following was shown:
\ms\nin
{\bf Theorem~6} (\cite[Theorem~2]{DSW}). {\it For $q=n-p>m=\bl[\nt\br]$, Conjecture~$2$ holds, if the following condition is satisfied$\,:$}
$$\h{\it The morphism $\beta_k^{(i)}$ in $(0.5)$ is surjective for $k=m(d-1)-p$ and $i=1$.}
\leqno(B')$$
\ms
Note that condition~$(B')$ is equivalent to condition~$(B)$ in {\it loc.~cit}. We show in this paper the following (see (2.3) below and also \cite{Di4}).
\ms\nin
{\bf Theorem~7.} {\it Condition~$(B')$ holds if $n$ is even or $n$ is odd and $q\ges m+[d/2]$.}
\ms
Combining Theorem~7 with Theorems~1 and 5, we get the following.
\ms\nin
{\bf Theorem~8.} {\it Conjecture~$2$ holds except for the case where $n$ is odd and $m<q<m+[d/2]$.}
\ms
The situation in the exceptional case is unclear (since condition~$(B')$ is only a sufficient condition), and Conjectures~1 and 2 are still open, see remarks after \cite[Theorem~2]{DSW}.
\sk
For the proof of Theorem~7, set
$${\rm def}_k\Si_f:=\dim{\rm Coker}\bl(\beta_k^{(1)}:R_k\to\mopl_{y\in\Sg\,Y}\,\OO_{X,y}/\mm_{X,y}\br),$$
so that condition~$(B')$ is equivalent to
$${\rm def}_{m(d-1)-p}\Si_f=0.
\leqno(B'')$$
It follows from the last isomorphisms of (0.2) that
$${\rm def}_k\Si_f=\tauY-\dim M''_{k+n+1}\q\h{with}\q\tauY:=\#|\Sg\,Y|.$$
By \cite[Theorem~3.1]{Di3} we have moreover
$$\dim\sN_{nd-n-1-k}={\rm def}_k\Si_f.
\leqno(0.6)$$
This also follows from \cite[Corollary~2]{DS2} asserting
$$\dim M''_k+\dim\sN_{nd-k}=\tauY.
\leqno(0.7)$$
These are closely related with \cite{BrHe}, \cite{Ei}, \cite{EyMe}, \cite{Pe}, \cite{Se}, \cite{vSt}, \cite{vStWa}. (Note that $n$ in \cite{DS2} is $n+1$ in this paper, and the grading of $\sN_{\ssb}$ in this paper is {\it shifted by} $d$ compared with $N_{\ssb}$ in \cite{DS2}; more precisely, $\sN_k=N_{k+d}$.)
\sk
It follows from (0.6) that condition~$(B'')$ is equivalent to
$$\sN_{(n-m)d+m-q-1}=0,
\leqno(B''')$$
since $nd-n-1-(m(d-1)-p)=(n-m)d+m-q-1$.
\sk
Here we have a generalization of results in \cite{DiSt1}, \cite{Di4} (which has been conjectured in \cite[Conjecture 3.15]{DP}, see also \cite{Kl}) as follows.
\ms\nin
{\bf Theorem~9.} {\it Assume all the singularities of $Y$ are isolated and weighted homogeneous. Then $\sN_k=0$ for $\kd<\alt_Y$, where $\alt_Y$ is as in $(1.1.3)$ below.}
\ms
If condition $(\ODP)$ is satisfied, then $\alt_Y=\nt$ and Theorem~9 was shown by \cite{DiSt1} in the $n$ even case, and by \cite{Di4} in the $n$ odd case.
(It is known in these cases that the bound on $k$ given in Theorem~9 is sharp; that is, there are examples with $\sN_k\ne 0$ for any integer $k\ges nd/2$, see {\it loc.~cit.})
Theorem~7 then follows by calculating the condition
$$(n-m)d+m-q-1<nd/2,
\leqno(0.8)$$
see (2.3) below. Theorem~9 is shown by using a recent result from \cite{DS2} (see (1.2.10) below) together with the Thom-Sebastiani type theorems, see (2.1) and (2.2) below.
\sk
The second named author was partially supported by JSPS Kakenhi 24540039.
\sk
In Section~1 we review some basics of the Hodge and pole order filtrations and pole order spectral sequences.
In Section~2 we prove the main theorems after showing Theorem~(2.1).
\bs\bs
\centerline{\bf 1. Preliminaries}
\bs\nin
In this section we review some basics of the Hodge and pole order filtrations and pole order spectral sequences.
\ms\nin
{\bf 1.1.~Cohomology of projective hypersurface complements.}
Let $Y$ be a hypersurface in $X:=\PP^n$. Set $U=X\setminus Y$. By Grothendieck there are canonical isomorphisms
$$H^j(U,\C)=H^j\bl(\Gamma(X,\Om_X^{\ssb}(*Y))\br)\q(j\in\Z),$$
where $\Gamma(X,\Om_X^{\ssb}(*Y))$ is the complex of rational differential forms on $X$ whose poles are contained in $Y$. By \cite{mhm}, there is a canonical Hodge filtration $F$ on $\OO_X(*Y)$ underlying a mixed Hodge module $j_*\Q_{h,U}[n]$ where $j:U\into X$ is the inclusion and $\Q_{h,U}[n]$ denotes the pure Hodge module of weight $n$ whose underlying $\Q$-complex is $\Q_U[n]$. By \cite[Proposition~2.2]{DS1}, we have
$$F^pH^j(U,\C)=H^j\bl(\Gamma(X,\Om_X^{\ssb}\otimes_{\OO_X}F_{\ssb-p}\OO_X(*Y))\br)\q(p,j\in\Z).
\leqno(1.1.1)$$
Here $F$ on the left-hand side coincides with the Hodge filtration of the canonical mixed Hodge structure on $H^j(U,\C)$ in \cite{De}. (This can be reduced to the case of the complement of a divisor with normal crossings on a smooth projective variety easily.)
\sk
Let $P$ be the pole order filtration on $\OO_X(*Y)$ defined by
$$P_p\OO_X(*Y):=\begin{cases}0&\h{if}\qt p<0,\\ \OO_X((p+1)Y)&\h{if}\qt p\ges 0.\end{cases}$$
Then
$$F_p\OO_X(*Y)\subset P_p\OO_X(*Y),$$
$$F_p\OO_X(*Y)|_{X\setminus\Sg\,Y}=P_p\OO_X(*Y)|_{X\setminus\Sg\,Y}.$$
Let $h_y$ be a local defining holomorphic function of $Y$ at $y$, and $b_{h_y}(s)$ be the $b$-function of $h_y$ which is normalized as in \cite{DS1}, \cite{bfut} so that
$$\widetilde{b}_{h_y}(s):=b_{h_y}(s)/(s+1)\in\C[s].$$
Let $\alt_{Y,y}$ be the minimal root of $\widetilde{b}_{h_y}(-s)$. If $(Y,y)$ is an isolated singularity defined locally by a weighted homogenous polynomial $h_y$ of weights $w_1,\dots,w_n$ (that is, $h_y$ is a linear combination of monomials $y_1^{m_1}\cdots y_n^{m_n}$ with $\sum_iw_im_i=1)$, then
$$\alt_{Y,y}=\msum_{i=1}^n\,w_i.
\leqno(1.1.2)$$
Set
$$\alt_Y:=\min_{y\in\Sg\,Y}\alt_{Y,y}.
\leqno(1.1.3)$$
By \cite{bfut} we have
$$F_p\OO_X(*Y)=P_p\OO_X(*Y)\q\h{if}\qt p<[\alt_Y].
\leqno(1.1.4)$$
This implies
$$F^pH^j(U,\C)=P^pH^j(U,\C)\q\h{if}\qt p>j-[\alt_Y],
\leqno(1.1.5)$$
where the filtration $P$ on $H^j(U,\C)$ is induced by $P$ on $\OO_X(*Y)$ by using the image of the right-hand side of (1.1.1) with $F$ replaced by $P$.
\sk
By \cite[Theorem~2.2]{DSW} we then get for $q=n-p<[\alt_Y]$
$$\Gr_F^pH^n(U,\C)=\Gr_P^pH^n(U,\C)=(R/J)_{(q+1)d-n-1}=M_{(q+1)d},
\leqno(1.1.6)$$
where $(R/J)_k$, $M_k$ are as in the introduction (although $Y$ may have arbitrary singularities).
\ms\nin
{\bf 1.2.~Pole order spectral sequences.}
In the notation of the introduction, we have the algebraic microlocal Gauss-Manin complex
$$(\Ct_f^{\ssb},\,\ddd-\dd_t\,\df\wedge)\q\h{with}\q \Ct_f^j=\Om^j[\dd_t,\dd_t^{-1}],$$
see \cite{DS2}.
Here $f$ may be any homogeneous polynomial of degree $d$.
Its cohomology groups $H^j(\Ct_f^{\ssb})$ are called the Gauss-Manin systems.
These are graded $\C$-vector spaces (where $\deg\dd_t=-d$), and there are isomorphisms
$$\aligned&H^{j+1}(\Ct_f^{\ssb})_k=\Ht^j(f^{-1}(1),\C)_{\la}\q\h{for}\qf\la=\exp(-2\pi ik/d),\\&\h{and}\q\q\Ht^j(f^{-1}(1),\C)_1=\Ht^j(U,\C),\endaligned
\leqno(1.2.1)$$
where $\Ht^j(f^{-1}(1),\C)_{\la}$ denotes the $\la$-eigenspace of the reduced Milnor cohomology under the monodromy.
It is well-known (see for instance \cite[1.3]{BS}) that there is a local system $L_k$ of rank 1 on $U$ such that $L_0=\C_U$, and
$$H^j(f^{-1}(1),\C)_{\la}=\begin{cases}H^j(U,L_k)&\h{if}\qt\la=\exp(-2\pi ik/d),\\\,0&\h{if}\qt\la^d\ne 1.\end{cases}
\leqno(1.2.2)$$
\sk
We have the pole order filtration $P'$ defined by
$$P'_p\,\Ct_f^j=\mopl_{i\les j+p}\,\Om^j\,\dd_t^i.$$
Set $P^{\prime\,p}=P'_{-p}$. There are isomorphisms
$$\dd_t^a:P^{\prime\,p}(\Ct_f^{\ssb})_k\simto P^{\prime\,p-a}(\Ct_f^{\ssb})_{k-ad}\q(a,k,p\in\Z).
\leqno(1.2.3)$$
We have the algebraic microlocal pole order spectral sequence
$$_kE_1^{p,j-p}=H^j\Gr_{P'}^p(\Ct_f^{\ssb})_k\Longrightarrow H^j(\Ct_f^{\ssb})_k,
\leqno(1.2.4)$$
which is a spectral sequence of graded $\C$-vector spaces.
The associated filtration $P'$ on $H^{j+1}(\Ct_f^{\ssb})_k$ is identified with the pole order filtration $P$ on the reduced Milnor cohomology groups $\Ht^j(f^{-1}(1),\C)_{\la}$ defined in \cite[Chapter~6]{Di1} up to the shift by one (that is, $P^{\prime\,p+1}=P^p$) via the isomorphism (1.2.1) for $k\in[1,d]$. This follows from \cite[Chapter~6, Theorem~2.9]{Di1} (see also \cite[Section~1.8]{DS1} for the case $j=n$). We thus get
$$P^{\prime p+1}H^{j+1}(\Ct_f^{\ssb})_k=P^p\Ht^j(f^{-1}(1),\C)_{\la}\q\h{for}\qf\la=\exp(-2\pi ik/d),\,k\in[1,d].
\leqno(1.2.5)$$
(For $k\notin[1,d]$, we have to use (1.2.3).)
\sk
By (1.1.5) we have moreover
$$F^p=P^p\q\h{on}\q\Ht^j(f^{-1}(1),\C)_1=\Ht^j(U,\C)\q\h{if}\qf p>j-[\alt_Y],
\leqno(1.2.6)$$
where $\alt_Y$ is as in (1.1.3). By definition we have
$$_kE_1^{n,0}=H^n\Gr_{P'}^n(\Ct_f^{\ssb})_k=\sN_k,\q_kE_1^{n+1,0}=H^{n+1}\Gr_{P'}^{n+1}(\Ct_f^{\ssb})_k=M_k,$$
where $\sN_k$, $M_k$ are as in the introduction. Set
$$\sN_k^{(r)}:={}_kE_{r}^{n,0},\q M_k^{(r)}:={}_kE_{r}^{n+1,0}\q(r\in[1,\infty]).
\leqno(1.2.7)$$
The differential $\ddd_r$ of the pole order spectral sequence is then identified by using (1.2.3) with the graded morphism
$$\ddd^{(r)}:\sN^{(r)}\to M^{(r)}\q(r\in[1,\infty)),$$
so that it preserves the grading up to a shift by $(r-1)d$, that is,
$$\ddd^{(r)}(\sN_k^{(r)})\subset M_{k-(r-1)d}^{(r)},
\leqno(1.2.8)$$
and moreover there are canonical isomorphisms
$$\sN^{(r+1)}={\rm Ker}\,\ddd^{(r)},\q M^{(r+1)}={\rm Coker}\,\ddd^{(r)}\q(r\in[1,\infty)).
\leqno(1.2.9)$$
\sk
If all the singularities of $Y$ are isolated and weighted homogeneous, then we have by \cite[Theorem~5.3]{DS2}
$$\sN^{(2)}_k={\rm Ker}(\ddd^{(1)}:\sN_k\to M_k)=0\q\h{if}\qf\kd<\alt_Y.
\leqno(1.2.10)$$
\sk
By (1.2.3), (1.2.5) we have in general
$$M_k^{(\infty)}=\Gr_P^pH^n(f^{-1}(1),\C)_{\la}\q\bl(p=\bl[n+1-\kd\br],\,\,\la=\exp(-2\pi ik/d)\br).
\leqno(1.2.11)$$
Indeed, this follows from (1.2.5) and the definition of the filtration $P'$ on $\Ct_f^{\ssb}$ given just before (1.2.3) if $k\in[1,d]$ so that $p=\bl[n+1-\kd\br]=n$. In the general case we also use the isomorphism (1.2.3) where $a\in\Z$ is chosen so that $k-ad\in[1,d]$.
\sk
By (1.2.6) we then get
$$M_{(q+1)d}^{(\infty)}=\Gr_F^pH^n(U,\C)\q\bl(q=n-p<[\alt_Y]\br).
\leqno(1.2.12)$$
\ms\nin
{\bf 1.3.~Thom-Sebastiani type theorems} (\cite[Section 4.9]{DS2}).
Let $h=f+g$ with $g$ a homogeneous polynomial of degree $d$ in variables $z_1,\dots,z_r$. In the notation of (1.2), there is a canonical isomorphism
$$(\Ct^{\ssb}_{h},P')=(\Ct^{\ssb}_f,P')\otimes_{\C[\dd_t,\,\dd_t^{-1}]}(\Ct^{\ssb}_g,P').
\leqno(1.3.1)$$
If $g$ has an isolated singularity at the origin, then
$$H^jGr^{P'}_k\Ct^{\ssb}_{g}=0\q(j\ne r,\,\,k\in\Z),
\leqno(1.3.2)$$
and we have the filtered quasi-isomorphisms
$$(\Ct^{\ssb}_{g},P')\simto H^r(\Ct^{\ssb}_{g},P')[-r],
\leqno(1.3.3)$$
$$(\Ct^{\ssb}_{h},P')\simto(\Ct^{\ssb}_{f},P')\otimes_{\C[\dd_t,\,\dd_t^{-1}]} H^r(\Ct^{\ssb}_{g},P')[-r].
\leqno(1.3.4)$$
\sk
By (1.3.4) the pole order spectral sequence for $h$ is isomorphic to a finite direct sum of shifted pole order spectral sequences for $f$ by choosing graded free generators of $H^r(\Ct^{\ssb}_{g},P')$ over $\C[\dd_t,\,\dd_t^{-1}]$. Here shifted means that the degrees of complex and filtration are shifted.
\sk
The $E_0$-complex of the spectral sequence is the direct sum of the graded quotients of the filtration $P'$, and is isomorphic to an infinite direct sum of the Koszul complexes $(\Om^{\ssb},\df\wedge)$ with grading shifted properly. We have the Thom-Sebastiani type theorem also for the Koszul complexes, see also \cite[Proposition~2.2]{DS2}.
\ms\nin
{\bf 1.4.~Spectrum.} We have the Steenbrink spectrum $\Sp(f)=\msum_{\al>0}\,n_{f,\al}\,t^{\al}$ defined by
$$\aligned n_{f,\al}:=\msum_j\,(-1)^j\dim\Gr^p_F\Ht^{n-j}(f^{-1}(1),\C)_{\la}\\
\h{with}\q p=[n+1-\al],\,\,\la=\exp(-2\pi i\al),\endaligned$$
where $\Ht^{n-j}(f^{-1}(1),\C)$ is the reduced cohomology, and $\Ht^{n-j}(f^{-1}(1),\C)_{\la}$ denotes the $\la$-eigenspace under the monodromy, see \cite{St2}, \cite{St3}.
We define the pole order spectrum $\Sp_P(f)=\msum_{\al>0}\,\Pn_{f,\al}\,t^{\al}$ by replacing $F$ with $P$.
\sk
There are refinements of the Steenbrink and pole order spectrum defined by
$$\Sp^j(f)=\msum_{\al>0}\,n^j_{f,\al}\,t^{\al},\q\Sp^j_P(f)=\msum_{\al>0}\,\Pn^j_{f,\al}\,t^{\al},$$
where
$$\aligned&n^j_{f,\al}:=\dim\Gr^p_F\Ht^{n-j}(f^{-1}(1),\C)_{\la}\\
\h{with}\q&p=[n+1-\al],\,\,\la=\exp(-2\pi i\al),\endaligned
\leqno(1.4.1)$$
and similarly for $\Pn^j_{f,\al}$ with $F$ replaced by $P$. By definition we have
$$\Pn^0_{f,k/d}=\dim M^{(\infty)}_k,\q \Pn^1_{f,k/d}=\dim\sN^{(\infty)}_{k-d}\,\br(=\dim N^{(\infty)}_k\br),
\leqno(1.4.2)$$
where $M^{(\infty)}_k$, $\sN^{(\infty)}_k$ are as in (1.2.7).
\sk
Since $F^p\subset P^p$ (see (1.7.4) below), we have
$$\Sp^j(f)-\Sp_P^j(f)=\msum_{\al,q}\,m^j_{\al,q}\,\bl(t^{\al}-t^{\al-q}\br),
\leqno(1.4.3)$$
where
$$\aligned&m^j_{\al,q}:=\dim\Gr_F^p\Gr_P^{p+q}\Ht^{n-j}(f^{-1}(1),\C)_{\la}\\
\h{with}\q&p=[n+1-\al],\,\,\la=\exp(-2\pi i\al),\,\,q\in\Z_{>0}.\endaligned$$
We also have
$$n^j_{f,\al}\les \Pn^j_{f,\al}\q\h{if}\qf n^j_{f,\al-k}=\Pn^j_{f,\al-k}\qf\h{for any}\qf k\in\Z_{>0}.
\leqno(1.4.4)$$
\sk
Let $h:=f+g$ with $g$ as in (1.3). Here we assume that $g$ has an isolated singularity at 0 so that
$$\Sp(g)=\Sp^0(g)=\Sp_P(g)=\Sp_P^0(g).$$
Then the Thom-Sebastiani type theorems imply
$$\Sp^j(h)=\Sp^j(f)\,\Sp^0(g),\q\Sp_P^j(h)=\Sp_P^j(f)\,\Sp^0(g),
\leqno(1.4.5)$$
where the product is taken in $\Q[t^{1/e}]$ for some positive integer $e$.
Indeed, the assertion for $\Sp^j(h)$ follows from \cite[Theorem 2]{MSS} (see also \cite{ts} for a different proof and \cite{SkSt} for the isolated singularity case). The assertion for $\Sp_P^j(h)$ follows from (1.3) by using (1.4.2).
\ms\nin
{\bf 1.5.~Spectral sequence.} There is a spectral sequence of mixed Hodge structures
$$_WE_1^{-i,i+j}=H^ji_0^*\Gr^W_i(\varphi_f\Q_{h,\C^{n+1}}[n])\Longrightarrow H^ji_0^*(\varphi_f\Q_{h,\C^{n+1}}[n]),
\leqno(1.5.1)$$
where $i_0:\{0\}\into\C^{n+1}$ denotes the inclusion, $\varphi_f\Q_{h,\C^{n+1}}[n]$ is a mixed Hodge module whose underlying $\Q$-complex is $\varphi_f\Q_{\C^{n+1}}[n]$, and $W$ is the weight filtration of the mixed Hodge module, see \cite{mhm}. Moreover there is a canonical isomorphism
$$H^ji_0^*(\varphi_f\Q_{\C^{n+1}}[n])=\Ht^{j+n}(f^{-1}(1),\Q),$$
compatible with the mixed Hodge structure.
\sk
Under the assumption~$(\ODP)$ in the introduction, we have the strict support decomposition of mixed Hodge modules
$$\Gr^W_i(\varphi_f\Q_{h,\C^{n+1}}[n])=M_{\{0\},i}\oplus\mopl_{y\in\Sg\,Y}M_{C(y),i},
\leqno(1.5.2)$$
where $C(y)\subset\C^{n+1}$ is the cone of $y\in\PP^n$, and $M_{Z,i}$ is a pure Hodge module of weight $i$ with strict support $Z=\{0\}$ or $C(y)$. Moreover the stalk of $M_{C(y),i}$ at any point of $C(y)\setminus\{0\}\cong\C^*$ is $\Q(-m)$ (up to a shift of complex by 1) if $i=2m+1$, and it vanishes otherwise, where $m:=[n/2]$ as in the introduction. (Indeed, these are well known if one restricts to the complement of the origin, and this implies the desired assertion by using the semisimplicity of pure Hodge modules.) The monodromy of $M_{C(y),2m+1}$ around the origin is the multiplication by $(-1)^{nd}$. Here we use a well-known relation with the Milnor monodromy on $(M_{C(y),2m+1})_z$ for $z\in C(y)\setminus\{0\}$, which is the multiplication by $(-1)^n$ by the assumption~$(\ODP)$. These imply
$$H^{-1}i_0^*\Gr^W_i(\varphi_f\Q_{h,\C^{n+1}}[n])=\begin{cases}\Q^{\,\tauY}(-m)&\h{if}\qt i=2m+1\qt\h{with}\qt nd\qt\h{even,}\\ \,0&\h{otherwise.}\end{cases}
\leqno(1.5.3)$$
$$H^0i_0^*\Gr^W_i(\varphi_f\Q_{h,\C^{n+1}}[n])=i_0^*M_{\{0\},i}.
\leqno(1.5.4)$$
Here $M_{\{0\},i}=(i_0)_*(i_0^*M_{\{0\},i})$, and $i_0^*M_{\{0\},i}$ can be identified with a pure Hodge structure $H_i$ of weight $i$.
\sk
We then get under the assumption $(\ODP)$
$$\h{$_WE_1^{-i,i+j}$ has weight $i+j$ for any $i,j$.}
\leqno(1.5.5)$$
and hence
$$\h{the spectral sequence (1.5.1) degenerates at $E_2$.}
\leqno(1.5.6)$$
Moreover the monodromical property of the weight filtration $W$ implies a decomposition of mixed Hodge structures
$$H_i=H_{i,1}\oplus H_{i,\ne1}$$
such that
$$N^k:H_{n+1+k,1}\simto H_{n+1-k,1}(-k),$$
$$N^k:H_{n+k,\ne 1}\simto H_{n-k,\ne 1}(-k),$$
as in the case of isolated singularities \cite{St2}, where $N=\log T_u$ with $T_u$ the unipotent part of the monodromy. By using (1.5.3--4) together with an argument similar to \cite{St2}, we then get
$$n_{f,\al}=n_{f,n+1-\al}\qf\h{for}\qf\al\ne\nt,\nt+1,
\leqno(1.5.7)$$
$$n_{f,n/2+1}=n_{f,n/2}-\tauY\qf\h{in the}\qt nd\qt\h{even case}.
\leqno(1.5.8)$$
For a similar assertion in a different setting, see \cite[Proposition 4.1]{Di2}. Note that the action of $N$ on $H^j(f^{-1}(1),\Q)$ is trivial (since $f$ is a homogeneous polynomial), and this implies a certain condition on the weight filtration $W$ on $\varphi_f\Q_{h,\C^{n+1}}[n]$.
By using the spectral sequence (1.5.1), it seems possible to give, for instance, another proof of \cite[Theorem~1.5]{DL} (where the equivariant Hodge-Deligne polynomial is essentially equivalent to the spectral pairs).
\ms\nin
{\bf 1.6.~Remark about (1.5.7--8).} These can be generalized to the case the hypothesis $(\ODP)$ is replaced by the condition that $Y$ has only hypersurface isolated singularities with {\it semi-simple} Milnor monodromies (for instance, a local defining function $h_y$ of $(Y,y)$ is weighted homogeneous for $y\in\Sg\,Y$). Indeed, let $\Sp(h_y)=\msum_{\al}\,n_{h_y,\al}\,t^{\al}$ be the spectrum of $h_y$, and define
$$\aligned&[\Sp(h_y)]_{(d)}:=\msum_{d\al\in\Z}\,n_{h_y,\al}\,t^{\al},\\&\widetilde{\rm Sp}(f):=\Sp(f)+\msum_{y\in\Sg\,Y}\,t\,[\Sp(h_y)]_{(d)}.\endaligned$$
Then $\widetilde{\rm Sp}(f)$ has symmetry with center $(n+1)/2$, that is, it is invariant by the automorphism of $\Z[t^{1/d},t^{-1/d}]$ defined by $Q(t)\mapsto t^{n+1}Q(t^{-1})$.
\ms\nin
{\bf 1.7.~Remark about the inclusion $F\subset P$.} Let $C_f^{\ssb}\subset\Ct_f^{\ssb}$ be the usual Gauss-Manin complex where $C_f^j=\Om^j[\dd_t]$ and the differential is $\ddd-\dd_t\,\df\wedge$. It has the pole order filtration $P'$ as in (1.2), and the inclusion induces the isomorphisms
$$P^{\prime\,p}H^{n+1}(C_f^{\ssb})_k\simto P^{\prime\,p}H^{n+1}(\Ct_f^{\ssb})_k\q(p\les n+1,k\in[1,d]),
\leqno(1.7.1)$$
since $\dd_t$ has degree $-d$.
\sk
The Gauss-Manin complex $C_f^{\ssb}$ can be identified with $f^{\D}_{*}\OO_{\X}$ up to the shift of complex by $n+1$, where $\X=\C^{n+1}$ and $f^{\D}_{*}$ denotes the direct image as an algebraic $\D$-module.
The Hodge filtration $F$ on the Gauss-Manin system $\Hc^0f^{\D}_{*}\OO_{\X}$ is defined by taking a relative compactification $\f:\XX\to S:=\C$ and using the factorization $f=\f\ssc j$ together with the Hodge filtration $F$ on $j^{\D}_*\OO_{\X}=j_*\OO_{\X}$, where $j:\X\into\XX$ is a compactification such that $D:=\XX\setminus\X$ is a divisor and $j_*$ is the direct image as an algebraic quasi-coherent sheaf. We have trivially the inclusions
$$F_p(j_*\OO_{\X})\subset j_*(F_p\OO_{\X}),
\leqno(1.7.2)$$
where $F_p\OO_{\X}=\OO_{\X}$ or 0.
The pole order filtration $P'$ on $H^{n+1}(C_f^{\ssb})$ is identified with the filtration induced by the direct image as a filtered $\D$-module of $j_*\OO_{\X}$ endowed with the filtration defined by the right-hand side of (1.7.2).
It is well-known (see for instance \cite[4.2.1]{DS2}) that $V$-filtration of Kashiwara and Malgrange on $H^{n+1}(C_f^{\ssb})$ is given by
$$V^{\al}H^{n+1}(C_f^{\ssb})=\mopl_{k\ges d\al}\,H^{n+1}(C_f^{\ssb})_k.
\leqno(1.7.3)$$
Combined with (1.2.1), these imply the inclusions
$$F^pH^n(f^{-1}(1),\C)_{\la}\subset P^pH^n(f^{-1}(1),\C)_{\la}.
\leqno(1.7.4)$$
Here we have the shift by one for the Hodge filtration $F$ as in (1.2.5).
This comes essentially from the transformation between filtered left and right $\D$-modules on $S$.
\ms\nin
{\bf 1.8.~Complement to the proof of \cite[Theorem~2.2]{DSW}.}
It does not seem to be necessarily easy to follow it, since it was written far too concisely. We give here an additional explanation as follows:
\sk
Using the vanishing theorem \cite[(2.1.1)]{DSW} together with the {\it strictness} of the Hodge filtration \cite{De}, we can express $\Gr^p_FH^n(U,\C)$ as the cokernel of the following morphism:
$$\ddd:\frac{\Gamma\bl(X,F_{n-p-1}\OO_X(*Y)\otimes_{\OO_X}\Om_X^{n-1}\br)}{\Gamma\bl(X,F_{n-p-2}\OO_X(*Y)\otimes_{\OO_X}\Om_X^{n-1}\br)}\to\frac{\Gamma\bl(X,F_{n-p}\OO_X(*Y)\otimes_{\OO_X}\Om_X^n\br)}{\Gamma\bl(X,F_{n-p-1}\OO_X(*Y)\otimes_{\OO_X}\Om_X^n\br)}
\leqno(1.8.1)$$
This can be computed further by using \cite[(2.1.2--5)]{DSW}. For $q:=n-p\les m$, it is identified, by using the morphism $\iota_{\xi}$ in {\it loc.~cit.}, with the cokernel of the morphism
$$\aligned\df\wedge:(\Om^n/f\Om^n)_{qd}\to(\Om^{n+1}/f\Om^{n+1})_{(q+1)d}\q\h{if}\qt q<m,\\
\df\wedge:(\Om^n/f\Om^n)_{qd}\to(I\Om^{n+1}/f\Om^{n+1})_{(q+1)d}\q\h{if}\qt q=m.\endaligned
\leqno(1.8.2)$$
Then \cite[Theorem~2.2]{DSW} follows. Here it seems important to write down the formulas (1.8.1--2) explicitly in order to understand the proof of \cite[Theorem~2.2]{DSW} properly.
\ms\nin
{\bf 1.9. Remark about \cite[Proposition 2.2]{DiSt2}.} It is shown there that $\Sp^1(f)=\Sp_P^1(f)$ in the case $Y$ has only isolated singularities.
The argument is closely related to \cite[Remark 4.4]{wh}, and is not quite trivial, since we do not know yet whether the restriction morphism is always {\it strictly compatible\,} with the pole order filtration in general. Indeed, we have a quite difficult problem that the {\it componentwise} strictness does not necessarily imply the strictness on the cohomology groups. More precisely, assume there is a morphism of filtered complexes
$$(C^{\ssb},P)\to(C'{}^{\ssb},P),$$
which is componentwisely strict, that is, the morphisms $(C^j,P)\to(C'{}^j,P)$ are strict ($\forall\,j$). However, the induced filtered morphism
$$\bl(H^j(C^{\ssb}),P\br)\to\bl(H^j(C'{}^{\ssb}),P\br)$$
is not necessarily {\it strict} (even in the case $C'{}^{j+1}=0$).
\sk
For instance, assume the morphism of complexes is given by
$$\begin{array}{cccccccccccc}
0&\to&\Q\oplus\Q&\buildrel{\alpha}\over\longrightarrow&\Q&\to&0\\ &&\,\,\,\downarrow{\!\scriptstyle\beta}&&\downarrow\\ 0&\to&\Q&\longrightarrow&0\end{array}
\leqno(1.9.1)$$
with $\Q\oplus\Q$ put at degree $0$. Assume $\Q\oplus 0$, ${\rm Ker}\,\alpha$, ${\rm Ker}\,\beta$ are different 1-dimensional subspaces. In particular, we have the isomorphism
$$H^0(C^{\ssb})\simto H^0(C'{}^{\ssb}).
\leqno(1.9.2)$$
Define the filtration $P^1$ by the morphism of subcomplexes
$$\begin{array}{cccccccccccc}
0&\to&\Q\oplus 0&\simto&\Q&\to&0\\ &&\,\,\,\downarrow{\!\scriptstyle\cong}&&\downarrow\\ 0&\to&\Q&\longrightarrow&0\end{array}
\leqno(1.9.3)$$
and assume $\Gr_P^jC^{\ssb}=\Gr_P^jC'{}^{\ssb}=0$ for $j\ne 0,1$, that is, $P^2C^{\ssb}=0$, $P^0C^{\ssb}=C^{\ssb}$, and similarly for $C'{}^{\ssb}$. We then get
$$\Gr_P^0H^0(C^{\ssb})=H^0(C^{\ssb})=\Q,\q\Gr_P^1H^0(C'{}^{\ssb})=H^0(C'{}^{\ssb})=\Q.
\leqno(1.9.4)$$
This implies, however, that $\bl(H^0(C^{\ssb}),P\br)\to\bl(H^0(C'{}^{\ssb}),P\br)$ is {\it nonstrict\,} by (1.9.2).

\bs\bs
\centerline{\bf 2. Proofs of the main theorems}
\bs\nin
In this section we prove the main theorems after showing Theorem~(2.1).
\ms\nin
{\bf 2.1.~Theorem.} {\it In the notation of $(1.2)$ we have
$$M_k^{(\infty)}=M_k\q\bl(\kd\les\alt_Y\br),
\leqno(2.1.1)$$
and hence}
$$\ddd^{(1)}:\sN_k\to M_k\qt\h{vanishes for}\qt\kd\les\alt_Y.
\leqno(2.1.2)$$
\ms\nin
{\it Proof.} It is enough to show (2.1.1). We prove it by using (1.1.6), (1.2.12) and the Thom-Sebastiani type theorems for the Gauss-Manin systems and the Koszul complexes as in (1.3). Set
$$h:=f+g\q\h{on}\q\C^{n+1}\times\C^r\q\h{with}\q g=\msum_{j=1}^r\,z_j^d\,$$
where
$$r=1\q\h{or}\q r=k_0:=\min\bl\{k\in\N\mid\alt_Y+\kd\ges[\alt_Y]+1\br\}.$$
Define $M_{h,k}$, $M_{h,k}^{(\infty)}$ by replacing $f$ with $h$, and similarly for $M_{g,k}=M_{g,k}^{(\infty)}$.
\sk
By the Thom-Sebastian type theorems as in (1.3), we have for $k\in\Z$
$$M_{h,k}=\mopl_{j\in\Z}\,(M_{f,k-j}\otimes_{\C}M_{g,j}),\q M_{h,k}^{(\infty)}=\mopl_{j\in\Z}\,(M_{f,k-j}^{(\infty)}\otimes_{\C}M_{g,j}),$$
where we denote $M_k$, $M_k^{(\infty)}$ by $M_{f,k}$, $M_{f,k}^{(\infty)}$ to avoid the ambiguity. We thus get
$$\dim M_{h,k}-\dim M_{h,k}^{(\infty)}=\msum_j\,(\dim M_{f,k-j}-\dim M_{f,k-j}^{(\infty)})\dim M_{g,j},
\leqno(2.1.4)$$
\sk
By (1.1.6) and (1.2.12) we have
$$M_{f,qd}=M_{f,qd}^{(\infty)}\q\h{for any integer}\qt q\les[\alt_Y].
\leqno(2.1.5)$$
Similarly we have
$$M_{h,qd}=M_{h,qd}^{(\infty)}\q\h{for any integers}\,\,\begin{cases}q\les[\alt_Y]&\h{if $\,r=1$},\\q\les[\alt_Y]+1&\h{if $\,r=k_0$}.\end{cases}
\leqno(2.1.6)$$
Indeed, if we set $Y':=\{h=0\}\subset\PP^{n+r}$, then it follows from (1.1.2) that
$$\alt_{Y'}=\begin{cases}\alt_Y+\tfrac{1}{d}\ges[\alt_Y]&\h{if}\qt r=1.\\ \alt_Y+\tfrac{r}{d}\ges[\alt_Y]+1&\h{if}\qt r=k_0,\end{cases}$$
by using a Thom-Sebastiani type theorem for $b$-functions, see \cite{micro}.
\sk
By definition we have
$$\dim M_{f,k}-\dim M_{f,k}^{(\infty)}\ges 0\q\h{for any}\qt k.$$
Moreover
$$\dim M_{g,j}>0\q\h{for}\qt j\in[r,r(d-1)],$$
since it is well-known that
$$\msum_j\,(\dim M_{g,j})\,t^j=(t+\cdots+t^{d-1})^r.
\leqno(2.1.7)$$
(This can be reduced to the case $r=1$ by using the Thom-Sebastiani type theorem.)
\sk
Using these non-negativity and strict positivity together with (2.1.4--6), we then get the equalities
$$\dim M_{f,k}-\dim M_{f,k}^{(\infty)}=0\q\h{for}\qt\kd\les\alt_Y.
\leqno(2.1.8)$$
Indeed, we first show the assertion for $\kd\les[\alt_Y]$ by using (2.1.5) and (2.1.6) for $r=1$. We then apply (2.1.6) for $r=k_0$ in case $\alt_Y\ne[\alt_Y]$.
\sk
We thus get a partial degeneration of the spectral sequence, and (2.1.1) follows.
This finishes the proof of Theorem~(2.1).
\ms\nin
{\bf 2.2.~Proof of Theorem~9.} This follows from (2.1.2) in Theorem~(2.1) together with \cite[Theorem 5.3]{DS2} (see (1.2.10)).
\ms\nin
{\bf 2.3.~Proof of Theorem~7.}
In the notation of the introduction, we see that condition~(0.8) becomes the following condition:
$$\begin{array}{ll}q>m-1\,\,&\h{if}\qt n=2m,\\q>m+d/2-1\,\,&\h{if}\qt n=2m+1.\end{array}$$
So Theorem~7 follows.
\ms\nin
{\bf 2.4.~Proof of Theorem~2.} By Theorem~(2.1) together with (1.2.8) and (1.2.11), we get the $E_2$-degeneration of the pole order spectral sequence together with the isomorphisms in (0.3) with $F$ replaced by $P$. The coincidence of $F$ and $P$ can be shown inductively by an argument similar to the proof of Theorem~(2.1) using the Thom-Sebastiani type theorems as in (1.4.5) together with (1.4.4). This finishes the proof of Theorem~2.
\ms\nin
{\bf 2.5.~Proof of Theorem~3.} Let $\al=\kd$ with $k\in\N$. The assertion is then equivalent to
$$\Pn_{f,\al}=\begin{cases}\ga_k&\h{if}\qf\al\les\nt,\\ \ga_k-\dim\sN_k&\h{if}\qf\nt<\al\les\nt+1,\\ \ga_k-(\dim\sN_k-\dim\sN_{k-d})&\h{if}\qf\al>\nt+1.\end{cases}$$
By the $E_2$-degeneration of the pole order spectral sequence in Theorem~2, we have
$$\Pn_{f,k/d}=\dim M_k^{(2)}-\dim\sN_{k-d}^{(2)}\,\bl(=\dim M_k^{(2)}-\dim N_k^{(2)}\br).
\leqno(2.5.1)$$
By \cite[Theorem~5.3]{DS2} (see (1.2.10)) we get
$$\dim M_k^{(2)}=\dim M_k-\dim\sN_k,\qf\dim\sN_k^{(2)}=0,\q\h{if}\qf \kd\ne\nt,
\leqno(2.5.2)$$
Here we state (0.4) again for the convenience of the reader (unless he has a very good memory since this is often used in this subsection):
$$\dim M_k=\ga_k+\dim\sN_{k-d}\q(k\in\Z).
\leqno(2.5.3)$$
We then easily see that (2.5.1), (2.5.2), and (2.5.3) imply the assertion for $\al\ne\nt,\,\nt+1$ by using Theorem~9 (that is, $\sN_k=0$ for $k<nd/2$).
\sk
Assume
$$k_0:=nd/2\in\N.$$
By (2.1.1), (2.5.3) and Theorem~9 we easily see that
$$\dim M_{k_0}^{(2)}=\dim M_{k_0}=\ga_{k_0},\q\dim\sN_{k_0-d}^{(2)}=0.
\leqno(2.5.4)$$
So the assertion for $\al=\nt$ follows.
\sk
By (2.1.2) we have
$$\sN_{k_0}^{(2)}=\sN_{k_0}.
\leqno(2.5.5)$$
We then easily see that the assertion for $\al=\nt+1$ follows by using (2.5.1), (2.5.2), and (2.5.3), where $k=k_0+d$ in order that $\al=\kd=\nt+1$. Indeed, there is a cancellation of $\dim\sN_{k-d}=\dim\sN_{k_0}$.
(Note that (2.5.5) together with the $E_2$-degeneration implies Proposition~1.)
This finishes the proof of Theorem~3.
\ms\nin
{\bf 2.6.~Proof of Theorem~4.} Let $\al=\kd$ with $k\in\N$. The assertion is then equivalent to
$$n_{f,\al}=\begin{cases}\ga_k&\h{if}\qf\al\les\nt\qf\h{or}\qf\al>\nt+1,\\ \ga_k-\tauY&\h{if}\qf\nt<\al\les\nt+1.\end{cases}$$
For $\al\les\nt$, this follows from Theorem~2 together with (0.4) and Theorem~9.
Combined with (1.5.7--8), it implies the assertion for $\al>\nt+1$.
By using (1.2.2) (see also \cite[1.4.2]{BS}), the assertion for $\nt<\al\les\nt+1$ is then reduced to
$$\chi(U)=\chi(U')-(-1)^n\tauY,
\leqno(2.6.1)$$
where $U'=\PP^n\setminus Y'$ with $Y'\subset\PP^n$ a nonsingular hypersurface of degree $d$.
\sk
For the proof of (2.6.1), we have the following well-known formula
$$\chi(Y)=\chi(Y')-(-1)^{n-1}\tauY.$$
Indeed, this can be shown by using a deformation of $Y$ since $Y$ has only isolated singularities. So (2.6.1) follows. This finishes the proof of Theorem~4.
\ms\nin
{\bf 2.7.~Proof of Theorem~5.} The morphisms are defined as is explained in Theorem~5. Since the pole order filtration degenerates at $E_2$ by Theorem~2, and the action of $\dd_t^{n-p}$ in (1.2.3) is an isomorphisms, it is enough to show
$${\rm Im}(\ddd^{(1)}:\sN_k\to M_k)\cap M'_k=0\q(k\in\N).
\leqno(2.7.1)$$
But this follows from the improved version of \cite[Theorem~5.3]{DS2} explained in \cite[Remark 5.6(i)]{DS2} if $k\ne nd/2$, and from Proposition~1 if $k=nd/2$. (Indeed, Proposition~1 is essentially equivalent to the vanishing of $\ddd^{(1)}:\sN_k\to M_k$ for $k=nd/2$.) This finishes the proof of Theorem~5.
\ms\nin
{\bf 2.8.~Proof of Proposition~2.} The assertion is equivalent to the inequality
$$n^0_{f,k/d}\ges\dim M'_k=\dim M-\dim M''_k=\ga_k+\dim\sN_{k-d}-\dim M''_k,
\leqno(2.8.1)$$
together with the equality for $\kd\in\bl(\nt,\nt+1\br]$, where the last equalities of (2.8.1) follow from (0.1) and (0.4).
\sk
We have a symmetry of $\{\dim M'_k\}_k$ with center $k=d(n+1)/2$ by \cite[Corollary~1]{DS2} (where $n$ means $n+1$ in this paper), and a partial symmetry of $\{n^0_{f,k/d}\}_k$ by (1.5.7). So it is enough to consider the following two cases
$$(1):\kd\les\nt+\ot,\q(2):\kd=\nt+1.$$
By Theorem~9 and (0.7) we have
$$\dim\sN_{k-d}=0\q\bl(\kd<\nt+1\br),
\leqno(2.8.2)$$
$$\dim M''_k=\tauY\q\bl(\kd>\nt\br).
\leqno(2.8.3)$$
So the assertion follows from Theorem~4 in the case (1).
\sk
In the case (2) we have by Proposition~(2.9) below
$$n^1_{f,k/d}=\dim\sN_{nd/2}.
\leqno(2.8.4)$$
Since $n^0_{f,k/d}=n_{f,k/d}+n^1_{f,k/d}$ by definition, the assertion follows from Theorem~4 by using (2.8.3--4), where we have a cancellation of $\dim\sN_{nd/2}$. This finishes the proof of Proposition~2.
\ms\nin
{\bf 2.9.~Proposition.} {\it In the notation of $(1.4)$, we have
$$\Sp^1(f)=\Sp_P^1(f)=\begin{cases}(\dim{}\sN_{nd/2})\,t^{n/2+1}&\h{if}\qt nd\qt\h{is even,}\\ \,0&\h{if}\qt nd\qt\h{is odd,}\end{cases}$$
and}
$$\aligned\Sp^0(f)&=\Sp(f)+\Sp^1(f),\\ \Sp_P^0(f)&=\Sp_P(f)+\Sp_P^1(f).\endaligned$$
\ms\nin
{\it Proof.} Since the last assertion follows from the definition, it is enough to show the first assertion, and the latter follows from Proposition~1 together with (1.5.3) and (2.1.2). This finishes the proof of Proposition~(2.9).
\ms\nin
{\bf 2.10.~Remark.} Without assuming condition~$(\ODP)$, the isomorphism in (0.3) holds for $\kd\les\alt_Y$ with $\alt_Y$ as in (1.1.3). This is shown by an argument similar to the proof of Theorem~2 in (2.4). Similarly we can show the $E_2$-degeneration of the pole order spectral sequence if all the singularities of $Y$ are isolated and weighted homogeneous and if $\alt_Y>(n-1)/2$.
Indeed, $\sN_k^{(2)}=0$ for $\kd\ges(n+1)/2$ by \cite[Theorem~5.3]{DS2}. Hence the $E_2$-degeneration follows from Theorem~(2.1) together with (1.2.8). (This assertion is recently proved without assuming $\alt_Y>(n-1)/2$, see \cite{wh}.)
\sk
It is rather difficult to generalize Theorems~3 and 4 even to the simple singularity case. For instance, if every singular point of $Y$ has type $A_{k_i}$ with $k_i\in\Z_{>0}$, then
$$\alt_Y=\tfrac{n-1}{2}+\tfrac{1}{k+1}\q\h{with}\q k=\max\{k_i\},$$
and it is not quite easy to determine $n_{f,\al}$ for $\al\in[\alt_Y,n-\alt_Y]\cup[\alt_Y+1,n+1-\alt_Y]$ if $k>1$.
\ms\nin
{\bf 2.11.~Examples} (see \cite[Examples~5.7]{DS2}). In the notation of (1.2.7) set
$$\mu_k=\dim M_k,\q{}^s\nu_k=\dim\sN_k,\q\mu^{(2)}_k=\dim M^{(2)}_k,\q{}^s\nu^{(2)}_k=\dim\sN^{(2)}_k.$$
\sk\nin
(i) $f=x^2y^2+x^2z^2+y^2z^2\qt$(three $A_1$ singularities in $\PP^2$) $\,n=2,\,d=4.$
$$\begin{array}{cccccccccccccccccc}
k\, &1 &2 &3 &4 &5 &6 &7 &8 &9 &10 &\cdots&\raise-3mm\h{ }\\
\ga_k & & &1 &3 &6 &7 &6 &3 &1\\
\mu_k & & &1 &3 &6 &7 &6 &3 &3 &3 &\cdots\\
{}^s\nu_k & & & & &2 &3 &3 &3 &3 &3 &\cdots \\
\mu^{\scriptscriptstyle(2)}_k & & &1 &3 &4 &4 &3\\
{}^s\nu^{\scriptscriptstyle\,(2)}_k \\
\Sp_P & & &1 &3 &4 &4 &3 & &\\
\Sp & & &1 &3 &3 &4 &3 & &1\\
\end{array}$$
\ms\nin
(ii) $f=xyz(x+y+z)\qt$(six $A_1$ singularities in $\PP^2$) $\,n=2,\,d=4.$
$$\begin{array}{cccccccccccccccccc}
k\, &1 &2 &3 &4 &5 &6 &7 &8 &9 &10 &\cdots&\raise-3mm\h{ }\\
\ga_k & & &1 &3 &6 &7 &6 &3 &1\\
\mu_k & & &1 &3 &6 &7 &6 &6 &6 &6 &\cdots\\
{}^s\nu_k & & & &3 &5 &6 &6 &6 &6 &6 &\cdots \\
\mu^{\scriptscriptstyle(2)}_k & & &1 &3 &1 &1\\
{}^s\nu^{\scriptscriptstyle\,(2)}_k & & & &3\\
\Sp_P & & &1 &3 &1 &1 & &-3 &\\
\Sp & & &1 &3 & &1 & &-3 &1\\
\end{array}$$

\end{document}